\documentclass[12pt]{article}
\usepackage{amsfonts,amsmath,amssymb,amscd}
\usepackage[russian]{babel}
\input epsf
\long\def\comment#1\endcomment{}

\begin{document}


\newpage
\centerline{\uppercase{\bf Простое доказательство теоремы Абеля}}
\smallskip
\centerline{\uppercase{\bf о неразрешимости уравнений в радикалах}
\footnote{ Заметка основана на занятиях кружка `Олимпиады и математика' и
выездных школ команды Москвы на Всероссийскую олимпиаду (http://www.mccme.ru/circles/oim/mat.htm).
Благодарю М. Вялого, А. Канеля-Белова, Г. Челнокова
и участников моих занятий за полезные обсуждения, а также М. Вялого за подготовку рисунков.}
}
\bigskip
\centerline{\bf А. Скопенков
\footnote{Поддержан грантом фонда Саймонса.
Московский Физико-Технический Институт,
Независимый Московский Университет,
Инфо: www.mccme.ru/\~{ }skopenko}
}

\bigskip
{\bf Теорема Абеля.}
{\it Калькулятор имеет кнопки
$$1,\quad +,\quad -,\quad \times,\quad :\quad\mbox{и}\quad \sqrt[n]{}\quad
\mbox{для любого $n$}.$$
Он оперирует с комплексными числами и при нажатии кнопки $\sqrt[n]{}$
выдает все значения корня (и как-то их случайно нумерует).
Калькулятор вычисляет числа с абсолютной точностью и имеет неограниченную память.
При делении на 0 он выдает ошибку и прекращает работу.

Тогда ни при каком $n\ge5$ не существует программы для этого калькулятора,
которая по коэффициентам многочлена $n$-й степени выдает конечное
множество чисел, содержащее все его корни.}

\smallskip
Комментарии по поводу понятия программы и другая формулировка приведены в пункте
'пут\"евые перестановки для программы и коммутаторы'.

В этой заметке мы покажем, как можно было бы придумать простое доказательство теоремы Абеля.
Мы следуем изложению [FТ], которое упрощено по сравнению с [A]
(не используется понятия римановой поверхности  и гомотопии).
В отличие от [FТ], здесь излагается способ {\it придумать} доказательство.
\footnote{
В одном месте (задача 7b, аналогично решаются задачи 8b и 9b) вместе с эвристическим рассуждением из [FT] приводится и более строгое (которое мне сообщил А. Канель со ссылкой на А. Ногина).
В отличие от [FТ], здесь используется понятие осторожного пути вместо рассуждений о перестановках башни значений радикальной формулы при обходе параметром замкнутого пути, не проходящего через особые точки радикальной формулы.
}
Для понимания самого доказательства достаточно прочитать определение пут\"евой перестановки из следующего пункта, пункт `план простого доказательства теоремы Абеля' и решения задач, на которые есть ссылки в этом пункте.

Заметим, что теорема (Галуа) о неразрешимости в радикалах {\it одного конкетного} уравнения [K, P, S2, T] --- более сильная.

Приводимое доказательство не используют терминов 'группа Галуа'
и 'расширение поля'.
Несмотря на отсутствие этих {\it терминов}, {\it идеи} приводимого
доказательства являются {\it отправными} для теории Галуа (более подробно см. [Ch, KS, S1, S3]).
Материал преподносится в виде задач, к которым даются указания.
(И отсутствие терминов, и присутствие задач, характерно не только для дзенских
монастырей, но и для серьезного изучения математики.)
В конце приводятся задачи для исследования.
Если условие задачи является формулировкой утверждения, то это утверждение
и надо доказать.

\smallskip
{\bf 1.}
(a) Теорему Абеля достаточно доказать для уравнений 5-й степени.

(b) Теорему Абеля достаточно доказать для уравнений $z^5-z+a=0$.

(На самом деле, и необходимо, поскольку любое уравнение пятой степени сводится к
написанному при помощи некоторой программы для нашего калькулятора.)

(c) (это не задача, а {\it загадка} [VINH])
Сформулируйте вещественный аналог теоремы Абеля.
Вытекает ли он из комплексного? А комплексный из вещественного?

\smallskip
{\bf 2.} (a) У калькулятора оторвали кнопки $\sqrt[n]{}$.
Теперь не существует программы для решения квадратного уравнения.

(b)* Не существует программы для решения кубического
уравнения, использующей извлечение корня только один раз.

(c)* Не существует программы для решения уравнения 4-й степени, использующей
извлечение корня не более двух раз.

\smallskip
Вряд ли у вас получится решить задачу 2c без знания дальнейшего материала.
К ней нужно возвращаться по мере его изучения.

\bigskip
{\bf Пут\"евые перестановки. }

В этом и следующих двух пунктах калькулятор не используется (поэтому даже читатель, у которого возникли вопросы о работе калькулятора, может смело решать задачи).

\smallskip
{\bf 3.} (a) Руководитель кружка двигался по единичной окружности на
комплексной плоскости, сделал один оборот и вернулся в исходную точку.
Он велел ученику двигаться на комплексной плоскости так, чтобы координата
$z$ ученика  в любой момент равнялась бы квадрату координаты $a$ руководителя.
Как пришлось двигаться ученику?

(b) На следующее занятие кружка на комплексную плоскость пришло два
ученика.
Руководитель встал в точке 1, а учеников расставил в точки 1 и $-1$.
Потом он велел каждому ученику двигаться так, чтобы координата $a$
руководителя в любой момент равнялась бы квадрату координаты $z$ ученика.
А сам пошел по единичной окружности на комплексной плоскости против
часовой стрелки, сделал один оборот и вернулся в исходную точку 1.
Как пришлось двигаться ученикам?
Где они оказались в конце занятия?

\smallskip
{\bf 4.} На следующее занятие кружка на комплексную плоскость пришло уже $n$
учеников.
Руководитель не растерялся, встал в точке 1, а учеников расставил
в точки $\cos(2\pi k/n)+i\sin(2\pi k/n)$, $k=1,2,\dots,n$.
Потом он сказал:
`Mоя координата равна $n$-й степени координаты любого из вас!
Двигайтесь так, чтобы сохранить это замечательное свойство.'

(a) Руководитель сам пошел по единичной окружности на комплексной плоскости против
часовой стрелки, сделал один оборот и вернулся в исходную точку 1.
Как пришлось двигаться ученикам?
Где они оказались в конце занятия?

(b) Руководитель называется {\it добрым}, если он не проходит через 0.
Для любого замкнутого маршрута доброго руководителя в конце занятия ученики
переставятся.

(c) Для произвольного замкнутого маршрута доброго руководителя если в конце ученик 1 оказался в точке $\cos(2\pi k/n)+i\sin(2\pi k/n)$, то ученик $\cos(2\pi/n)+i\sin(2\pi/n)$
оказался в точке $\cos(2\pi(k+1)/n)+i\sin(2\pi(k+1)/n)$.

(d) Для произвольного замкнутого маршрута доброго руководителя
в конце занятия ученики переставятся по степени некоторого цикла.


\smallskip
\centerline{\epsfbox{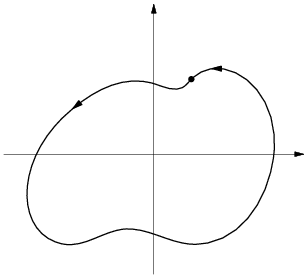}\qquad\qquad \epsfbox{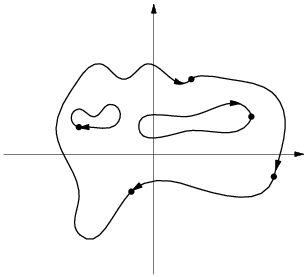}}
\centerline{изменение параметра $a$ \qquad\qquad соотвествующее изменение корней}
\centerline{Рисунок 1}
\smallskip

\smallskip
{\it Определение пут\"евой перестановки.}
Пусть $p_a(z)$ --- семейство многочленов степени $n$, коэффициенты которого (но не степень) непрерывно зависят от параметра $a$.
Пусть уравнение $p_{a_0}(z)=0$ имеет $n$ различных корней, которые обозначены $z_1,\dots,z_n$.
Будем изменять параметр $a$ (руководитель) вдоль некоторого непрерывного замкнутого
пути с началом и концом в $a_0$.
Пусть для любой точки $a$ этого пути уравнение $p_a(z)=0$ имеет $n$ различных корней.
Будем двигать $i$-го ученика, начиная в $z_i$, и так,
чтобы в каждый момент времени его координата была корнем уравнения
 Тогда в конце движения ученики переставятся.
\footnote{Для первого знакомства с приводимыми идеями читателю полезно
воспользоваться без доказательства существованием такого движения учеников (ср. задачи 4bcd),
а также использовать аналогичные наглядные соображения при решении задач 7b, 8b и 9b ниже.
Строгое обоснование  вытекает из {\it комплексной теоремы о неявной функции}. }
Назовем полученную перестановку $n$-элементного множества {\bf пут\"евой} для семейства уравнений
$p_a(z)=0$ и данной нумерации корней уравнения $p_{a_0}(z)=0$.

\smallskip
Это действительно перестановка, т.е. два ученика
не могут в конце оказаться в одной точке (задача 4b).
Цикл $(12\dots n)$ является путевым для $p_a(z)=z^n-a$ и `естественной' нумерации его корней (задача 4a); все пут\"евые перестановки являются степенями этого цикла (задача 4d).

\smallskip
{\bf 5.} (a) Тождественная перестановка пут\"евая для любого семейства многочленов и любой  нумерации корней.

(b) Как перестановка, отвечающая вышеописанному замкнутому пути, зависит от начальной нумерации корней?

\smallskip
С этого места мы пропускаем слова `для некоторой нумерации корней', говоря о
пут\"евых перестановках.

\smallskip
{\bf Ответы и указания. }

\smallskip
{\bf 4.} (c) Если $x(t)$ --- закон движения первого ученика, то
$x(t)(\cos(2\pi /n)+i\sin(2\pi /n))$ --- закон движения второго.

(d) Следует из (c).

\smallskip
{\bf 5.} (a) Если параметр $a$ в процессе движения стоит на месте, то каждый корень остается на месте.

(b) Ответ: при перенумерации, заданной перестановкой $\alpha$, перестановка $\sigma$, отвечающая вышеописанному замкнутому пути, перейдет в перестановку $\alpha\sigma\alpha^{-1}$.

\bigskip
{\bf Зачем нужны пут\"евые перестановки?}

Если бы удалось доказать, что любая пут\"евая перестановка для семейства
уравнений, разрешимого в радикалах, циклическая, то для доказательства теоремы
Абеля достаточно было бы привести пример семейства уравнений и нециклической
пут\"евой для него перестановки.
Однако перестановка $(13)(24)=(1234)^2$ пут\"евая для $p_a(z)=z^4-a$.

Хорошо было бы найти другое свойство пут\"вых перестановок для
уравнений, разрешимых в радикалах, которое не выполняется для пут\"евых перестановок произвольных уравнений.
Этого сделать не получится, ибо
{\it любая перестановка является пут\"евой
для некоторого семейства уравнений, разрешимого в радикалах, и некоторой нумерации корней.}
(Действительно, если перестановка является произведением $k$ циклов длин $n_1,\dots,n_k$,
то искомое семейство уравнений --- $\Pi_{s=1}^k((z-n_s)^{n_s}-a)$.)

И все-таки мы докажем теорему Абеля.
Мы найдем свойство {\it множества} всех пут\"евых перестановок для уравнений, разрешимых в радикалах, не выполненное для {\it множества} всех пут\"евых перестановок произвольных уравнений.

Покажем отправную идею на примере решения задачи 2a.
(Это решение сложнее придуманного вами ранее, но зато обобщается до доказательства теоремы Абеля.)

\smallskip
{\bf 6.} (a) Для каких $a$ уравнение $z^2-2z+a=0$ имеет ровно два корня?
(Здесь и далее имеются в виду комплексные корни.)

{\it Ответ.} Для $a\ne1$.

(b) Как переставляются корни уравнения $z^2-2z+a=0$ при следующем изменении параметра $a$?

Сначала от 0 до 0.99 по отрезку,

потом по окружности радиуса 0.01, обходящей вокруг точки 1 один раз

наконец, обратно от 0.99 до 0 по отрезку.

(с) Если $q(a)$ является отношением многочленов от $a$, то только тождественная перестановка является пут\"евой для $p_a(z)=z-q(a)$ (здесь путь параметра $a$ не проходит через нули числителя и знаменателя).

(d) Выведите из (b,c) решение задачи 2a.

\smallskip
{\bf Ответы и указания. }

\smallskip
Обозначим $e^{it}:=\cos t+i\sin t$. (В настоящем тексте это нужно воспринимать именно как обозначение. Свойство $e^{i(t_1+t_2)}=e^{it_1}e^{it_2}$ следует из формул для синуса и косинуса суммы.)

\smallskip
{\bf 6.} (b)
Сначала корни приближаются к $1$ с разных сторон.
Чтобы корень двигался по закону $z(t)=1+\varepsilon e^{it}$, параметр $a$
должен двигаться по закону $a(t)=2z(t)-z^2(t)=1-\varepsilon^2 e^{2it}$.
Поэтому возьмем $\varepsilon=0.1$.
Когда $a$ совершит один оборот, каждый из корней совершит
пол-оборота, т.е. они поменяются местами.
Далее корни снова удалятся от 1.

{\it Ответ:} корни меняются местами.

\bigskip
{\bf Пут\"евые перестановки для уравнений 3-й, 4-й и 5-й степени.}

\smallskip
{\bf 7.} (a) Для каких $a$ уравнение $z^3-3z+a=0$ имеет ровно три корня?

{\it Указание.} Проще всего решать эту задачу при помощи следующих фактов:

$\bullet$ любой многочлен степени $n$ имеет ровно $n$ комплексных корней, и

$\bullet$ любой кратный корень многочлена является также корнем его
производной.

{\it Ответ.} Для $a\ne\pm2$.

(b) Как переставляются корни уравнения $z^3-3z+a=0$ при
следующем изменении параметра $a$?

Сначала от $0$ до $2-\delta_1$ по отрезку,

потом по кривой, `обходящей вокруг точки $2$ один раз',

наконец, обратно от $2-\delta_2$ до 0 по отрезку.

({\it В задачах 7bc, 8bc и 9bc выберите сами кривую и малые числа
$\delta_1,\delta_2$, чтобы было удобно находить перестановку.})

(с) Как переставляются корни
уравнения $z^3-3z+a=0$ при
следующем изменении параметра $a$?

Сначала от $0$ до $-2-\delta_1$ по отрезку,

потом по кривой, `обходящей вокруг точки $-2$ один раз',

наконец, обратно от $-2-\delta_2$ до 0 по отрезку.

(d) Для $p_a(z)=z^3-3z+a$ все перестановки путевые.

\smallskip
{\bf 8.} (a) Для каких $a$ уравнение $z^4-4z+a=0$ имеет ровно четыре корня?

Ответ. Для $a\ne3,3\alpha,3\alpha^2$, где $\alpha=(1+i\sqrt3)/2$.

(b) Как переставляются корни
уравнения $z^4-4z+a=0$ при следующем изменении параметра $a$?

Сначала от $0$ до $3-\delta_1$ по отрезку,

потом по кривой, `обходящей вокруг точки $3$ один раз',

наконец, обратно от $3-\delta_2$ до 0 по отрезку.

(с) Как переставляются корни
уравнения $z^4-4z+a=0$ при следующем изменении параметра $a$?

Сначала от $0$ до $(3-\delta_1)\alpha$ по отрезку,

потом по кривой, `обходящей вокруг точки $3\alpha$ один раз',

наконец, обратно от $(3-\delta_2)\alpha$ до 0 по отрезку.

(d) Для $p_a(z)=z^4-4z+a$ все перестановки пут\"евые.

\smallskip
{\bf 9.} (a) Для каких $a$ уравнение $z^5-5z+a=0$ имеет ровно пять корней?

Ответ. Для $a\ne4,4i,-4,-4i$.

(b) Как переставляются корни
уравнения $z^5-5z+a=0$ при следующем изменении параметра $a$?

Сначала от $0$ до $4-\delta_1$ по отрезку,

потом по кривой, `обходящей вокруг точки $4$ один раз',

наконец, обратно от $4-\delta_2$ до 0 по отрезку.

(с) Как переставляются корни
уравнения $z^5-5z+a=0$ при следующем изменении параметра $a$?

Сначала от $0$ до $(4-\delta_1)i$ по отрезку,

потом по кривой, `обходящей вокруг точки $4i$ один раз',

наконец, обратно от $(4-\delta_2)i$ до 0 по отрезку.

(d) Для $p_a(z)=z^5-5z+a$ все перестановки путевые.

\smallskip
{\bf 10.} Найдите все пут\"евые перестановки для

(a) $p_a(z)=z^4+2(1-2a)z^2+1$. \quad (b) $p_a(z)=(z^3-a)^3-a(a-1)$.

\smallskip
{\bf Ответы и указания. }

\smallskip
{\bf 7.} (b) Ответ: корни $0,\sqrt3$ поменяются местами, корень $-\sqrt3$ остается на месте.

Зададим движение корня 0, по нему восстановим движение параметра $a$, а затем
движение остальных корней.
См. рис. 2.

\smallskip
\centerline{\epsfbox{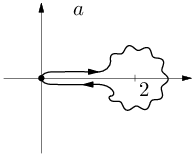}\qquad \qquad\epsfbox{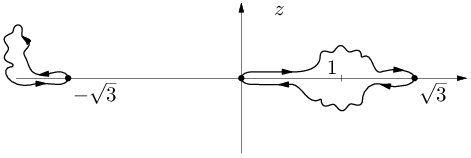}}
\centerline{изменение параметра $a$ \qquad\qquad соотвествующее изменение корней}
\centerline{Рисунок 2}
\smallskip

Сначала первый корень $z=0$ приближается к $1$ слева и превращается в $1-\delta_1$.
При этом $a=3z-z^3$ приближается к $2$ слева.
Значит, второй корень $\sqrt3$ приближается к $1$ справа и превращается в $1+\delta_2$.
А третий корень $-\sqrt3$ остается отрицательным.
См. рис. 3.

\smallskip
\centerline{\epsfbox{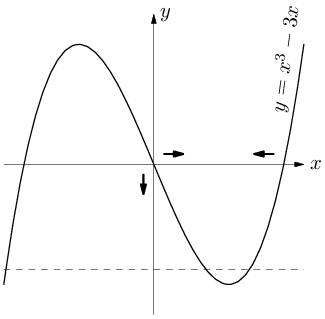}}
\centerline{Рисунок 3}
\smallskip


Затем пусть первый корень идет в $1+\delta_2$ по некоторой кривой, близкой
к $1$ и не пересекающей вещественной оси нигде, кроме своих концов.
Тогда второй корень остается близким к 1, а третий корень остается отрицательным.
Значит, корень второй корень придет в $1-\delta_1$.
\footnote{Вот неформальная иллюстрация этого движения.
Пусть первый корень движется по закону $z(t)=1+\delta_1e^{it}$.
Тогда параметр $a$ движется по закону
$a(t)=3z(t)-z^3(t)=2-3\delta_1^2 e^{2it}-\delta_1^3 e^{3it}$.
Так как $\delta_1$ мало, эта кривая близка к окружности $2-3\delta_1^2 e^{2it}$.
А значит, второй корень придет примерно в точку $1-\delta_1$ (ибо третий корень остается  далеко).
}

Далее первый корень приближается к $\sqrt3$ слева, примерно повторяя начальную часть движения второго корня в противоположном направлении.
Значит, второй корень приближается к 0 справа, примерно повторяя начальную часть движения первого корня в противоположном направлении.
А третий корень остается отрицательным.

(c) Указание. Используйте нечетность.

Ответ: корни $0,-\sqrt3$ поменяются местами, корень $\sqrt3$ остается на месте.

(d) Докажите, что любая перестановка 3-элементного множества представляется в виде
композиции транспозиций $(12)$ и $(13)$.

\smallskip
{\bf 8.} (b)
Аналогично задаче 7b корни $0,\sqrt[3]4$ поменяются местами.
Остальные два корня в процессе движения не пересекают ось $Ox$ и поэтому
в конце движения будут на своих прежних местах.

Ответ: корни $0,\sqrt[3]4$ поменяются местами, остальные два корня остаются
на месте.

(c) $p_{\alpha a}(\alpha z)=\alpha p_a(z)$.

Ответ: корни $0,\sqrt[3]4\alpha$ поменяются местами, остальные два корня
остаются на месте.

(d) Докажите, что любая перестановка 4-элементного множества представляется в виде
композиции транспозиций $(12)$, $(13)$ и $(14)$.

\smallskip
{\bf 9.} (b) Аналогично задачам 7b и 8b.
Корни $0,\sqrt[4]5$ поменяются местами, а
остальные три корня вернутся на свои прежние места.
(Ибо корни $i\sqrt[4]5$ и $-i\sqrt[4]5$ в процессе движения не пересекают
вещественной оси, а корень $-\sqrt[4]5$ остается вещественным отрицательным.)

(d) Докажите, что любая перестановка 5-элементного множества представляется в виде
композиции транспозиций $(12)$, $(13)$, $(14)$ и $(15)$.
Для этого докажите, что

$\bullet$ любая перестановка 5-элементного множества
является композицией циклов,

$\bullet$ любой цикл является композицией транспозиций,

$\bullet$ любая транспозиция является композицией транспозиций $(12)$, $(13)$, $(14)$ и $(15)$.

\bigskip
{\bf Осторожные пути.}

\smallskip
{\bf 11.} (a) Пусть $q(a)$ --- отношение многочленов от $a$.
Докажите, что все пут\"евые перестановки для $p_a(z)=z^n-q(a)$ являются степенью
некоторого одного цикла.

(b)* Пусть существует программа для решения уравнения $p_a(z)=0$, использующая
извлечение корня только один раз.
Обязательно ли все пут\"евые перестановки являются степенью
некоторого одного цикла?

\smallskip
Наш калькулятор имеет неприятную особенность: результат вычислений не всегда однозначно определяется вводимыми данными (например, программа, выдающая {\it первое} значение $\sqrt1$, будет случайно выдавать 1 или $-1$).
Эту неприятность можно преодолеть, осознав, что теорему Абеля достаточно доказать для `симметричных' программ для нашего калькулятора (или, эквивалентно, для похожего калькулятора, оперирующего с {\it множествами} комплексных чисел).
Но все равно будет непросто дать определение пут\"евой перестановки {\it для программы}, которое необходимо для решения задачи 11b.
Мы поступим по-другому.


{\it Радикальной формулой относительно $a_0,\dots,a_n$} называется (упорядоченный) набор рациональных функций (т.е. отношений многочленов) $p_1,\dots,p_s$ от $n+1,\dots,n+s$ переменных, соответственно, и целых положительных чисел $k_1,\dots,k_s$.
Для радикальной формулы определим выражения
$$z_1,\dots,z_s\quad\text{формулой}\quad z_j^{k_j}=p_j(a_0,\dots,a_n,z_1,\dots,z_{j-1}),\quad  j=1,2,\dots,s.$$
Уравнение $a_nz^n+a_{n-1}z^{n-1}+\dots+a_1z+a_0=0$ с комплексными переменными коэффициентами
называется {\it разрешимым в радикалах}, если существует радикальная формула относительно $a_0,\dots,a_n$, для которой любой корень уравнения является одним из значений одного из выражений $z_1,\dots,z_s$.
Ср. [FT, Chapter 5].

Мы докажем теорему Абеля в следующей эквивалентной форме:
{\it ни при каком $n\ge5$ уравнение $a_nz^n+a_{n-1}z^{n-1}+\dots+a_1z+a_0=0$ с комплексными переменными коэффициентами не является разрешимым в радикалах.}

Для
доказательства теоремы Абеля полезно следующее понятие.
Назовем замкнутый путь на плоскости {\bf осторожным} для данной радикальной формулы,
в которой $a_0,\dots,a_n$ зависят от параметра $a$, если при изменении параметра $a$ вдоль этого пути каждое значение каждого выражения $z_s$ возвращается на место.

\smallskip
{\bf 12.} (a) Если путь является осторожным относительно каждой из двух радикальных формул
$p_1,\dots,p_s;\ k_1,\dots,k_s$ и $q_1,\dots,q_t;\ l_1,\dots,l_t$,
то он является осторожным относительно их {\it суммы}
$$p_1,\dots,p_s,\ q_1,\dots,q_t,\ p_s+q_t;\quad k_1,\dots,k_s,\ l_1,\dots,l_t,\ 1.$$
\quad
(b) Определите разность, произведение и частное радикальных формул.
Докажите для них аналог предыдущего пункта.

(c) {\it $n$-й степенью} замкнутого пути называется новый замкнутый путь, который получается прохождением исходного пути $n$ раз.
Если путь является осторожным относительно радикальной формулы
$p_1,\dots,p_s;\ k_1,\dots,k_s$, то его $n$-я степень является осторожной относительно
радикальной формулы $p_1,\dots,p_s,\ z_s;\ k_1,\dots,k_s,\ n$.

(d) $k_1\cdot\ldots\cdot k_s$-я степень любого пути является осторожной относительно радикальной формулы $p_1,\dots,p_s;\ k_1,\dots,k_s$.

(e) Если семейство уравнений $p_a(z)=0$ разрешимо при помощи радикальной формулы
$p_1,\dots,p_s;\ k_1,\dots,k_s$, то $k_1\cdot\ldots\cdot k_s$-я степень любой пут\"евой перестановки тождественна.

\smallskip
Последним утверждением нельзя воспользоваться для доказательства теоремы Абеля, ибо любая перестановка в некоторой степени равна тождественной.

{\it Коммутатором} двух замкнутых путей называется новый замкнутый путь, который получается последовательным прохождением

первого пути,

второго пути,

первого пути в обратную сторону,

второго пути в обратную сторону.

\smallskip
{\bf 13.} (a) Если оба пути --- осторожные относительно радикальной формулы
\linebreak
$p_1,\dots,p_s;\ k_1,\dots,k_s$, то их коммутатор --- осторожный относительно
радикальной формулы $p_1,\dots,p_s,\ z_s;\ k_1,\dots,k_s,\ n$.

(b) Если существует программа для решения уравнения $p_a(z)=0$,
использующая одноэтажные извлечения корней, то любые две пут\"евые
перестановки коммутируют (т.е. $\sigma\tau=\tau\sigma$).


(c) Не существует программы для решения кубического
уравнения, использующей одноэтажные извлечения корней.
(Сравните с задачей 2b и вашим решением ее.)

\smallskip
{\bf 14.}
(a) Какое условие на множество пут\"евых перестановок следует из
существования программы для решения уравнения $p_a(z)=0$,
использующей не более, чем двухэтажное извлечение корня?

(b) Выведите из вашего решения пункта (a) и задачи 8d решение задачи 2c.

\smallskip
{\bf 15.}
(a) Какое условие на множество пут\"евых перестановок следует из
существования программы для решения уравнения $p_a(z)=0$,
использующей не более, чем трехэтажное извлечение корня?

(b) Существуют две не коммутирующие перестановки 5-элементного множества,
каждая из которых является коммутатором некоторых двух произведений
коммутаторов.

(c) Выведите из (a,b) и задачи 9 необходимость
наличия хотя бы трех этажей в
программе, якобы решающей семейство уравнений $z^5-z+a=0$.

(d) Докажите теорему Абеля.

\smallskip
{\bf Ответы и указания. }

\smallskip
{\bf 11.} (a) Аналогично задачам 4сd.

\smallskip
{\bf 13.} (a) Значение выражения $z_s$ возвращается на место в результате обхода числом $a$ каждого из двух данных замкнутых путей $L_1$ и $L_2$.
Поэтому $n$ значений выражения $z_{s+1}$ имеют вид
$x,x\varepsilon,x\varepsilon^2,x\varepsilon^3,\dots,x\varepsilon^{n-1}$ для некоторого
$x$ и $\varepsilon:=\cos(2\pi /n)+i\sin(2\pi /n)$.
Аналогично задаче 4d для любого замкнутого пути $L$ найдется такое $k(L)$, что в результате изменения параметра $a$ вдоль этого пути число $x\varepsilon^s$ переходит в число $x\varepsilon^{s+k(L)}$.
Поэтому в результате прохождения параметром $a$ коммутатора путей $L_1$ и $L_2$ число $x\varepsilon^s$ переходит в число $x\varepsilon^{s+k(L_1)+k(L_2)-k(L_1)-k(L_2)}$.

(b) Следует из (a).

\smallskip
{\bf 14.} (a) {\it Подсказка.}
Условие $\sigma\tau=\tau\sigma$ (на перестановки) равносильно
тождественности перестановки $\sigma\tau\sigma^{-1}\tau^{-1}$.
Эта перестановка называется {\it коммутатором} перестановок $\sigma,\tau$.
Если имеется двухэтажная формула, то для пут\"евых перестановок $\sigma,\tau$
коммутатор (т.е. перестановка $\sigma\tau\sigma^{-1}\tau^{-1}$)
может не быть тождественной.
Однако для коммутаторов выполняется некоторое условие.
Найдите его!

{\it Ответ.}
Если существует программа для решения уравнения $p_a(z)=0$,
использующая не более, чем двухэтажное извлечение корня, то коммутаторы
пут\"евых перестановок коммутируют.
(Даже произведения коммутаторов пут\"евых перестановок коммутируют.)

{\it Доказательство} аналогично решению задачи 13a.

\smallskip
{\bf 15.}
(a) Если существует программа для решения уравнения $p_a(z)=0$,
использующая не более, чем трехэтажное извлечение корня, то коммутаторы коммутаторов
пут\"евых перестановок коммутируют.
(Даже произведения коммутаторов произведений коммутаторов пут\"евых перестановок коммутируют.)
Доказательство аналогично задачам 13ab и 14a.

(b) Пример можно придумать напрямую или доказав, что любая четная
перестановка является произведением коммутаторов.
(Определение четной перестановки напомнено ниже.)

\bigskip
{\bf План простого доказательства теоремы Абеля. }

Мы довольно долго {\it придумывали} доказательство теоремы Абеля.
{\it Изложить} же доказательство можно совсем коротко.
(Освобождение доказательства от деталей, возникших при его придумывании и не нужных для него самого --- важная часть его проверки.)
Приведем план такого изложения.

Теорема Абеля вытекает из следующих трех лемм.
Для их формулировки введем следующее определение (которое поможет нам коротко произносить громоздкие конструкции, возникшие в задачах 13b, 14a и 15a).
Для данного семейства уравнений $p_a(z)=0$ назовем {\it 0-пут\"евыми} пут\"евые перестановки.
Если уже определены $n$-пут\"евые перестановки, то назовем перестановку {\it $(n+1)$-пут\"евой},
если она представляется в виде композиции коммутаторов некоторых $n$-пут\"евых перестановок.

\smallskip
{\bf Лемма о коммутаторах.} {\it Если существует программа для решения семейства уравнений $p_a(z)=0$, использующая не более, чем $n$-кратное извлечение корня, то лишь тождественная перестановка является $n$-пут\"евой для этого семейства.}

\smallskip
{\bf Лемма о примере уравнения.}
{\it Существует такое семейство $p_a(z)$ уравнений 5-й степени, множество всех пут\"евых перестановок которого совпадает с множеством всех перестановок 5-элементного множества.}

\smallskip
{\bf Лемма о четных перестановках.} {\it Любая четная перестановка 5-элементного
множества является произведением коммутаторов четных перестановок.}

\smallskip
Напомним, что перестановка называется {\it четной},
если она представляется в виде произведения циклов длины 3.
(Это определение равносильно общепринятому.)

Лемма о коммутаторах следует из задачи 13a (аналогично задачам 13b, 14a, 15a).
Лемма о примере уравнения следует из задачи 9d.
Доказательство леммы о четных перестановках --- задача (указание: достаточно
доказать, что таковым является цикл длины 3).

\bigskip
{\bf Задачи для исследования. }

В математике имеется много результатов, непосредственно связанных с теоремой Абеля.
См., например, [Kh], [PS].

Следующие задачи показывают, что метод Феррари
для решения уравнения 4-й степени `самый простой', а формула дель
Ферро-Кардано-Тартальи в виде
$$x=\frac12(\sqrt[3]{-a+\sqrt{a^2-4}}-\sqrt[3]{a+\sqrt{a^2-4}})$$ для решения кубического уравнения $x^3-3x+a=0$ --- не `самая простая'.

\smallskip
{\bf 16.} Существует программa для калькулятора, строящая по числу $a$
конечное множество, содержащее все корни уравнения $x^3-3x+a=0$, и
содержащая только одно извлечения корня из выражения, содержащего корни.


\smallskip
{\bf 17.} Не существует

(a) формулы вида $z=\sqrt[k]{p+\sqrt[l]{q+\sqrt[m]r}}+\sqrt[n]{s+\sqrt[o]t}$ для
решения уравнения $x^4-4x+a=0$ ни для каких целых положительных $k,l,m,n,o$ и
рациональных функций $p,q,r,s,t$ от $a$.

(b) программы, строящей по числу $a$
конечное множество, содержащее все корни уравнения $x^4-4x+a=0$, и содержащей
только одно `трехэтажное' извлечения корня.

\smallskip
{\bf 18.} Существует ли программа для {\it вещественного} аналога
калькулятора, определенного в начале заметки, находящая все {\it вещественные}
корни

(a) уравнения $x^3+px+q=0$ по его коэффициентам $p,q$?

(b) уравнения $x^4+px^2+qx+r=0$ по его коэффициентам $p,q,r$?

(c) уравнения $x^5+px^3+qx^2+rx+s=0$ по его коэффициентам $p,q,r,s$?

\smallskip
{\bf 19.} Добавим к калькулятору кнопку, выдающую по числу $\cos\alpha$ все
значения числа $\cos(\alpha/5)$. Появится ли программа для решения уравнения
5-й степени?

\smallskip
Задачи 16 и 17 несложны.
(Задачу 16 можно решить, не используя изложенных выше идей, а задачу 17 --- вряд ли.)
К сожалению, автору не удалось найти в литературе ответы на естественные
вопросы задач 18 и 19 (хотя, видимо, ответы известны специалистам).

\smallskip
{\bf Ответы и указания. }

\smallskip
{\bf 16.} $b:=\frac12\sqrt[3]{-a+\sqrt{a^2-4}}$ и $x:=b-\frac1b$.

\smallskip
{\bf 17.} Если бы такая формула/программа существовала, то все коммутаторы произведений коммутаторов (перестановок 4-элементного множества) были бы степенью некоторой одной нестареповки.
А это не так.

\bigskip
{\bf Литература.}

[A] В.Б. Алексеев, Теорема Абеля. М: Наука, 1976.


[Ch] Г.Р. Челноков, Основы теории Галуа в интересных задачах,
\linebreak
http://www.mccme.ru/circles/oim/materials/grishalois.pdf.
(засмотрено 11.11.2010)

[FT] D. Fuchs, S. Tabachnikov,  Mathematical Omnibus. AMS, 2007.

[K] В.А. Колосов, Теоремы и задачи алгебры, теории чисел и комбинаторики.
М: Гелиос, 2001.

[Kh] А.Г. Хованский, Топологическая теория Галуа, Москва, МЦНМО, 200?

[KS] П. Козлов и А. Скопенков, В поисках утраченной алгебры: в направлении
Гаусса (подборка задач), Мат. Просвещение, 12 (2008), 127--144,
http://arxiv.org/abs/0804.4357 v1

[P] В.В. Прасолов, Многочлены, М: МЦНМО, 2003.
http://www.mccme.ru/prasolov

[PS] В.В. Прасолов и Ю.П. Соловьев,
Эллиптические функции и алгебраические уравнения. М.: Факториал, 1997.
http://www.mccme.ru/prasolov

[S1] А. Скопенков, Философски-методическое отступление, в кн. 	
Сборник материалов московских выездных математических школ.
Под редакцией А. Заславского, Д. Пермякова, А. Скопенкова, М. Скопенкова и
А. Шаповалова, Москва, МЦНМО, 2009.
\linebreak
http://www.mccme.ru/circles/oim/mvz.pdf (засмотрено 20.08.2010)

[S2] А. Скопенков, Еще несколько доказательств из Книги: разрешимость
и неразрешимость уравнений в радикалах, http://arxiv.org/abs/0804.4357 v2

[S3] A. Skopenkov, Yet another proof from the Book:
the Gauss theorem on regular polygons, http://arxiv.org/abs/0908.2029

[T] В.М. Тихомиров, Абель и его великая теорема, Квант, 2003, N1.

[VINH] Виро О.Я., Иванов О.А., Нецветаев Н.Ю., Харламов В.М.	
Элементарная топология, М: МЦНМО, 2010.

\end{document}